\RequirePackage{fix-cm} 
\documentclass[a4paper, twoside, reqno, dvips, 12pt]{amsart}
\usepackage{fixltx2e}   

\usepackage{etex}

\usepackage[latin1]{inputenc}
\usepackage[T1]{fontenc}

\usepackage{esint}
\usepackage{dsfont}
\usepackage{xspace}
\usepackage{amsgen}
\usepackage{amsthm}
\usepackage{amssymb}
\usepackage{amsmath}
\usepackage{wasysym}
\usepackage{upgreek}
\usepackage{amsfonts}
\usepackage{stmaryrd}
\usepackage{mathtools}
\usepackage{textgreek}

\usepackage{relsize}
\usepackage[mathcal, mathscr]{euscript}

\usepackage{mathrsfs}
\DeclareMathAlphabet{\mathscrbf}{OMS}{mdugm}{b}{n}

\usepackage{a4wide}

\usepackage[dvipsnames, table]{xcolor}
\definecolor{bckg}{RGB}{20.8, 20.8, 20.8}
\definecolor{oneblue}{rgb}{0.0, 0.0, 0.85}
\definecolor{Lightblue}{RGB}{214, 214, 214}
\definecolor{bluepigment}{rgb}{0.2, 0.2, 0.6}
\definecolor{charcoal}{rgb}{0.21, 0.27, 0.31}
\definecolor{denimblue}{rgb}{0.08, 0.38, 0.74}
\definecolor{Lightgray}{rgb}{0.89, 0.89, 0.89}
\definecolor{darkgrey}{rgb}{0.273, 0.281, 0.30}
\definecolor{darkelectricblue}{rgb}{0.33, 0.41, 0.47}

\usepackage[sort&compress, comma, square, numbers]{natbib}

\usepackage{psfrag}
\usepackage{graphicx}
\usepackage{subfigure}
\usepackage{morefloats}
\usepackage{indentfirst}

\usepackage{acronym}
\usepackage{microtype}
\usepackage[labelsep=period,%
            labelfont={bf,sf,color=bluepigment},%
            justification=raggedright]{caption}

\usepackage[perpage, symbol]{footmisc}

\usepackage[usenames, dvipsnames, pdf]{pstricks}
\usepackage{epsfig}
\usepackage{pst-grad} 
\usepackage{pst-plot} 

\usepackage[colorlinks,
           urlcolor=oneblue,
           linkcolor=denimblue,
           citecolor=NavyBlue,
           bookmarksopen=false,
           pdfpagemode=UseNone,
           pagebackref]{hyperref}

\usepackage[explicit]{titlesec}

\titleformat{\section}[block]
  {\color{NavyBlue}\Large\sffamily\bfseries}
  {}
  {0.0em}
  {\colorbox{bckg!5}{\strut\parbox{\dimexpr\linewidth-2\fboxsep\relax}{\thesection. #1}}}
  [\vspace*{0.33em}]

\titleformat{name=\section,numberless}[block]
  {\color{NavyBlue}\Large\sffamily\bfseries}
  {}
  {0.0em}
  {\colorbox{bckg!5}{\strut\parbox{\dimexpr\linewidth-2\fboxsep\relax}{#1}}}
  [\vspace*{0.33em}]

\titleformat{\subsection}
  {\color{NavyBlue}\large\sffamily\bfseries}
  {}
  {0.0em}
  {\colorbox{bckg!5}{\parbox{\dimexpr\linewidth-2\fboxsep\relax}{\thesubsection. #1}}}
  [\vspace*{0.33em}]

\titleformat{name=\subsection,numberless}
  {\color{NavyBlue}\Large\sffamily\bfseries}
  {}
  {0em}
  {\colorbox{bckg!5}{\parbox{\dimexpr\linewidth-2\fboxsep\relax}{#1}}}
  [\vspace*{0.33em}]

\titleformat{\subsubsection}
  {\color{bluepigment}\sffamily\normalsize\bfseries}
  {\thesubsubsection}
  {0.5em}
  {#1}
  [\vspace*{0.33em}]

\titleformat{\paragraph}[runin]
  {\color{bluepigment}\sffamily\small\bfseries}
  {}
  {0em}
  {#1}

\titlespacing{\section}{0.0em}{1.5em plus 2pt minus 2pt}%
{1.0em plus 2pt minus 2pt}[0em]
\titlespacing{\subsection}{0.5em}{1.5em plus 2pt minus 2pt}%
{1.0em}[0em]
\titlespacing{\subsubsection}{0.5em}{1.5em plus 2pt minus 2pt}%
{1.0em plus 2pt minus 2pt}[0em]

\usepackage{titletoc}

\contentsmargin{0.5em}
\newlength{\tocsep} 
\titlecontents{section}[\tocsep]
  {\addvspace{10pt}\bfseries\sffamily}
  {\contentslabel[\thecontentslabel]{\tocsep}}
  {}
  {\ \titlerule*[0.75pc]{.}\ \thecontentspage}
  []
\titlecontents{subsection}[\tocsep]
  {\addvspace{8pt}\sffamily}
  {\contentslabel[\thecontentslabel]{\tocsep}}
  {}
  {\ \titlerule*[0.5pc]{.}\ \thecontentspage}
  []
\titlecontents*{subsubsection}[\tocsep]
  {\addvspace{2pt}\footnotesize\sffamily}
  {}
  {}
  {\ \titlerule*[0.35pc]{.}\ \thecontentspage}
  [\\*]

\makeatletter
\def\@setauthors{%
  \begingroup
  \def\thanks{\protect\thanks@warning}%
  \trivlist
  \centering\footnotesize \@topsep30\p@\relax
  \advance\@topsep by -\baselineskip
  \item\relax
  \author@andify\authors
  \def\\{\protect\linebreak}%
  \textsc{\normalsize\textcolor{darkelectricblue}{\authors}}%
  \ifx\@empty\contribs
  \else
    ,\penalty-3 \space \@setcontribs
    \@closetoccontribs
  \fi
  \endtrivlist
  \endgroup
}
\def\@settitle{\begin{center}%
  \baselineskip14\p@\relax
    \bfseries
    \textsc{\Large\textcolor{charcoal}{\@title}}
  \end{center}%
}
\makeatother

\usepackage{enumitem}
\setlist[description]{%
  topsep=30pt,               
  itemsep=5pt,               
  font={\bfseries\sffamily\color{NavyBlue}}, 
}

\usepackage{fancyhdr}
\usepackage{lastpage}

\newcommand*\Title{\textcolor{bluepigment}{Numerical solution of Feller's equation}}
\newcommand*\Authors{\textcolor{bluepigment}{D.~Dutykh}}
\newcommand*{\plogo}{\textcolor{gray}{{\texttt{arXiv.org} / \textsc{hal}}}} 

\numberwithin{equation}{section}

\newtheorem{remark}{Remark}

\newcommand{\C}{\mathds{C}}
\newcommand{\N}{\mathds{N}}
\newcommand{\R}{\mathds{R}}

\renewcommand{\tau}{\uptau}
\newcommand{\ud}{\mathrm{d}}
\newcommand{\ue}{\mathrm{e}}
\newcommand{\ui}{\mathrm{i}}
\newcommand{\x}{\mathscr{X}}
\newcommand{\y}{\mathscr{Y}}
\newcommand{\F}{\mathfrak{F}}
\newcommand{\G}{\mathfrak{G}}
\renewcommand{\O}{\mathcal{O}}
\renewcommand{\P}{\mathscr{P}}
\renewcommand{\leq}{\leqslant}
\renewcommand{\geq}{\geqslant}
\newcommand{\X}{\text{\textxi}}
\newcommand{\eps}{\upvarepsilon}
\newcommand{\const}{\mathrm{const}}

\newcommand{\D}{\mathscr{D}}
\newcommand{\Co}{\mathscr{C}}

\DeclareMathOperator{\E}{E}
\DeclareMathOperator{\Ei}{Ei\,}

\DeclareMathOperator{\M}{\mathcal{M}\,}
\DeclareMathOperator{\U}{\mathcal{U}\,}

\newcommand{\ie}{\emph{i.e.}\xspace}
\newcommand{\eg}{\emph{e.g.}\xspace}
\newcommand{\etc}{\emph{etc.}\xspace}

\renewcommand{\sim}{\thicksim}

\newcommand{\abs}[1]{\lvert\, #1\, \rvert}

\newcommand{\norm}[1]{\lVert\, #1\, \rVert}

\newcommand{\pd}[2]{\frac{\partial\/ #1}{\partial\/ #2}}
\newcommand{\PD}[2]{\dfrac{\partial\/ #1}{\partial\/ #2}}
\newcommand{\od}[2]{\frac{\mathrm{d}\/ #1}{\mathrm{d}\/#2}}

\newcommand{\odd}[2]{\dfrac{\mathrm{d}\/ #1}{\mathrm{d}\/#2}}
\newcommand{\eqdef}{\mathop{\stackrel{\,\mathrm{def}}{:=}\,}}

\newcommand{\half}{{\textstyle{1\over2}}}

\usepackage{acronym}

\begin{document}

\title[\Title]{Numerical simulation of Feller's diffusion equation}

\author[D.~Dutykh]{Denys Dutykh$^*$}
\address{\textbf{D.~Dutykh:} Univ. Grenoble Alpes, Univ. Savoie Mont Blanc, CNRS, LAMA, 73000 Chamb\'ery, France and LAMA, UMR 5127 CNRS, Universit\'e Savoie Mont Blanc, Campus Scientifique, F-73376 Le Bourget-du-Lac Cedex, France}
\email{Denys.Dutykh@univ-savoie.fr}
\urladdr{http://www.denys-dutykh.com/}
\thanks{$^*$ Corresponding author}



\begin{titlepage}
\thispagestyle{empty} 
\noindent
{\Large Denys \textsc{Dutykh}}\\
{\it\textcolor{gray}{CNRS--LAMA, Universit\'e Savoie Mont Blanc, France}}
\\[0.08\textheight]

\vspace*{2.5cm}

\colorbox{Lightblue}{
  \parbox[t]{1.0\textwidth}{
    \centering\huge\sc
    \vspace*{0.7cm}
    
    \textcolor{bluepigment}{Numerical simulation of Feller's diffusion equation}
    
    \vspace*{0.7cm}
  }
}

\vfill 

\raggedleft     
{\large \plogo} 
\end{titlepage}


\newpage
\thispagestyle{empty} 
\par\vspace*{\fill}   
\begin{flushright} 
{\textcolor{denimblue}{\textsc{Last modified:}} \today}
\end{flushright}


\newpage
\maketitle
\thispagestyle{empty}


\begin{abstract}

This article is devoted to \textsc{Feller}'s diffusion equation which arises naturally in probabilities and physics (\eg wave turbulence theory). If discretized naively, this equation may represent serious numerical difficulties since the diffusion coefficient is practically unbounded and most of its solutions are weakly divergent at the origin. In order to overcome these difficulties we reformulate this equation using some ideas from the \textsc{Lagrangian} fluid mechanics. This allows us to obtain a numerical scheme with a rather generous stability condition. Finally, the algorithm admits an elegant implementation and the corresponding \textsc{Matlab} code is provided with this article under an open source license.


\bigskip
\noindent \textbf{\keywordsname:} \textsc{Feller} equation; parabolic equations; \textsc{Lagrangian} scheme; numerical simulation; probability density function \\


\end{abstract}


\newpage
\tableofcontents
\thispagestyle{empty}


\newpage
\section{Introduction}

The celebrated \textsc{Feller} equation was introduced in two seminal papers published by William \textsc{Feller} (1951/1952) in Annals of Mathematics \cite{Feller1951, Feller1952}. These publications studied mathematically (and, henceforth, gave the name) to the following equation\footnote{To be more accurate, W.~\textsc{Feller} studied the following equation \cite{Feller1951}: \begin{equation*}
  p_{\,t}\ =\ \bigl[\,a\,x\,u\,\bigr]_{\,x\,x}\ -\ \bigl[\,(c\ +\ b\,x)\,u\,\bigr]_{\,x}\,,
\end{equation*} where $a\ >\ 0$ and $0\ <\ x\ <\ +\,\infty\,$.}:
\begin{equation}\label{eq:feller}
  p_{\,t}\ +\ \F_{\,x}\ =\ 0\,, \qquad \F\,(p,\,x,\,t)\ \eqdef\ -x\cdot\bigl(\gamma\,\,p\ +\ \eta\,p_{\,x}\bigr)\,,
\end{equation}
where the subscripts $t\,$, $x$ denote the partial derivatives, \ie $(\cdot)_{\,t}\ \eqdef\ \PD{(\cdot)}{t}\,$, $(\cdot)_{\,x}\ \eqdef\ \PD{(\cdot)}{x}\,$. Two parameters $\gamma$ and $\eta\ >\ 0$ can be time dependent in some physical applications, even if in this study we assume they are constants, for the sake of simplicity\footnote{The numerical method we are going to propose can be straightforwardly generalized for this case when $\gamma\ =\ \gamma\,(t)$ and $\eta\ =\ \eta\,(t)\,$. Moreover, \textsc{Feller}'s processes with time-varying coefficients were studied recently in \cite{Masoliver2016}.}. Equation~\eqref{eq:feller} can be seen as the \textsc{Fokker}--\textsc{Planck} (or the forward \textsc{Kolmogorov}) equation with $\gamma\,x$ being the drift and $\eta\,x$ being the diffusion coefficients. One can notice also that Equation~\eqref{eq:feller} becomes singular at $x\ =\ 0$ and $x\ =\ +\,\infty\,$. We remind that practically important solutions to \textsc{Feller}'s equation might be unbounded near $x\ =\ 0\,$. In order to attempt at solving Equation~\eqref{eq:feller}, one has to prescribe an initial condition $p\,(x,\,0)\ =\ p_{\,0}\,(x)$ presumably with a boundary condition at $x\ =\ 0\,$. A popular choice is to prescribe the homogeneous boundary condition $p\,(0,\,t)\ \equiv\ 0\,$. For this choice of the boundary condition it is not difficult to show that the \textsc{Feller} equation dynamics would preserve solution positivity provided that $p_{\,0}\,(x)\ \geq\ 0$ (see Appendix~\ref{app:pos} for a proof). The solution norm is also preserved (see Appendix~\ref{app:norm}). Moreover, \textsc{Feller} using the \textsc{Laplace} transform techniques has shown in \cite{Feller1951} that the initial condition $p_{\,0}\,(x)$ determines uniquely the solution. In other words, \emph{no} boundary condition at $x\ =\ 0$ should be prescribed. This conclusion might appear, perhaps, to be counter-intuitive.

The great interest in \textsc{Feller}'s equation can be explained by its connection to \textsc{Feller}'s processes, which can be described by the following stochastic differential \textsc{Langevin} equation\footnote{The stochastic differential equations is understood in the sense of \textsc{It\={o}}.}:
\begin{equation*}
  \ud X_{\,t}\ =\ -\,\gamma\,X_{\,t}\,\ud t\ +\ \sqrt{2\,\eta\,X_{\,t}}\,\ud\mathcal{W}_{\,t}\,,
\end{equation*}
where $\mathcal{W}_{\,t}$ is the standard \textsc{Wiener} process, \ie $\xi\,(t)\ \eqdef\ \odd{\mathcal{W}_{\,t}}{t}$ is zero-mean \textsc{Gaussian} white noise, \ie
\begin{equation*}
  \langle \xi\,(t)\rangle\ =\ 0\,, \qquad
  \langle \xi\,(t)\,\xi\,(s)\rangle\ =\ \delta\,(t\ -\ s)\,,
\end{equation*}
where the brackets $\langle \cdot\rangle$ denote an ensemble averaging operator. Then, the Probability Density Function (PDF) $p\,(x,\,t;\,x_{\,0})$ of the process $X\,(t)\,$, \ie 
\begin{equation*}
  \mathbb{P}\,\bigl\{x\ <\ X\,(t)\ <\ x\ +\ \ud x\ \vert\ X\,(0)\ =\ x_{\,0}\bigr\}\ \equiv\ p\,(x,\,t;\,x_{\,0})\,\ud x
\end{equation*}
satisfies Equation~\eqref{eq:feller} with the following initial condition \cite{Gardiner2004}:
\begin{equation*}
  p_{\,0}\,(x)\ =\ \delta\,(x\ -\ x_{\,0})\,, \qquad x_{\,0}\ \in\ \R^{\,+}\,.
\end{equation*}
The point $x\ =\ 0$ is a singular boundary that the process $X\,(t)$ cannot cross. The \textsc{Feller} process is a continuous representation of branching and birth-death processes, which never attains negative values. This property makes it an ideal model not only in physical, but also in biological and social sciences \cite{Murray2007, Gan2015, Masoliver2016}.

As a general comprehensive reference on generalized \textsc{Feller}'s equations we can mention the book \cite{Lehnikg1993}. Since at least a couple of years there is again a growing interest for studying equation \eqref{eq:feller}. Some singular solutions to \textsc{Feller}'s equation with constant coefficients were constructed in \cite{Gan2015} via spectral decompositions. \textsc{Feller}'s equation and \textsc{Feller}'s processes with time-varying coefficients were studied analytically (always using the \textsc{Laplace} transform) and asymptotically in \cite{Masoliver2016}.

Recently, the \textsc{Feller} equation was derived in the context of the weakly interacting random waves dominated by four-wave interactions \cite{Choi2005}. Wave Turbulence\footnote{We could define the Wave Turbulence (WT) as a physical and mathematical study of systems where random and coherent waves co-exist and interact \cite{Zakharov1992}.} (WT) is a common name for such processes \cite{Zakharov1992}. In WT the \textsc{Feller} equation governs the PDF of squared \textsc{Fourier} wave amplitudes, \ie $x\ \sim\ \abs{a}^{\,2}\,$. In \cite{Choi2005} some steady solutions to this equation with finite flux in the amplitude space were constructed\footnote{There is a misprint in \cite[p.~366]{Choi2005}. To obtain mathematically correct solutions one has to define $n_{\,k}\ \eqdef\ \dfrac{\eta}{\gamma}$ on the line below Equation (14).}. See also \cite[Chapter~11]{Nazarenko2011} for a detailed discussion and interpretations.

The behaviour of solutions $p\,(x,\,t)$ for large $x$ describes the appearance probability of extreme waves. In the context of ocean waves, these extreme events are known as \emph{rogue} (or \emph{freak}) waves \cite{Dysthe2008}. In the WT literature, any noticeable deviation from the \textsc{Rayleigh} distribution for $x\ \gg\ 1$ is referred to as the \emph{anomalous probability distribution} of large amplitude waves \cite{Choi2005}. For \textsc{Gaussian} wave fields all statistical properties can be derived from the spectrum. However, the PDFs and other higher order moments are compulsory tools to study such deviations.

The present study focuses on the numerical discretization and simulation of \textsc{Feller} equation. The naive approach to solve this equation numerically encounters notorious difficulties. The first question, which arises is what is the (numerical) boundary condition to be imposed at $x\ =\ 0\,$? Moreover, one can notice that Equation \eqref{eq:feller} is posed on a semi-infinite domain. There are three main strategies to tackle this difficulty:
\begin{enumerate}
  \item Map $\R^{\,+}$ on a finite interval $[\,0,\,\ell\,]$
  \item Use spectral expansions on $\R^{\,+}$ (\eg \textsc{Laguerre} or associated \textsc{Laguerre} polynomials)
  \item Replace (truncate) $\R^{\,+}$ to $[\,0,\,L\,]\,$, with $L\ \gg\ 1\,$.
\end{enumerate}
In most studies the latter option is retained by imposing some appropriate boundary conditions at the artificial boundary $x\ =\ L\,$. In our study we shall propose a method, which is able to handle the semi-infinite domain $\R^{\,+}$ without any truncations or simplifications. Finally, the diffusion coefficient in the \textsc{Feller} equation \eqref{eq:feller} is unbounded. If the domain is truncated at $x\ =\ L\,$, then the diffusion coefficient takes the maximal value $\nu_{\,\max}\ \eqdef\ \eta\,L\ \gg\ 1\,$, which depends on the truncation limit $L$ and it can become very large in practice. We remind also that explicit schemes for diffusion equations are subject to the so-called \textsc{Courant}--\textsc{Friedrichs}--\textsc{Lewy} (CFL) stability conditions \cite{Courant1928}:
\begin{equation*}
  \Delta t\ \leq\ \frac{\Delta x^{\,2}}{2\,\nu_{\,\max}}\,.
\end{equation*}
Taking into account the fact that $\nu_{\,\max}$ can be arbitrarily large, no explicit scheme can be usable with \textsc{Feller} equation in practice. Moreover, the dynamics of the \textsc{Feller} equation spreads over the space $\R^{\,+}$ even localized initial conditions. In general, one can show that the support of $p\,(x,\,t)\,$, $t\ >\ 0$ is strictly larger\footnote{Using modern analytical techniques it is possible to show even sharper results on the solution support, see \eg \cite{Carrillo2004}.} than the one of $p\,(x,\,0)\,$. It is the so-called \emph{retention property}. Thus, longer simulation times require larger domains. For all these reasons, it becomes clear that numerical discretization of the \textsc{Feller} equation requires special care.

In this study we demonstrate how to overcome this assertion as well. The main idea behind our study is to bring together PDEs and Fluid Mechanics. First, we observe that the classical \textsc{Eulerian} description is not suitable for this equation, even if the problem is initially formulated in the \textsc{Eulerian} setting. Consequently, the \textsc{Feller} equation will be recast in special \emph{material} or the so-called \textsc{Lagrangian} variables\footnote{It is known that both \textsc{Eulerian} and \textsc{Lagrangian} descriptions were proposed by the same person --- Leonhard \textsc{Euler}.}, which make the resolution easier and naturally adaptive \cite[Chapter~7]{Gosse2013}.

The present manuscript is organized as follows. The symmetry analysis of Equation~\eqref{eq:feller} is performed in Section~\ref{sec:sym}. Then, the governing equation is reformulated in \textsc{Lagrangian} variables in Section~\ref{sec:reform}. The numerical results are presented in Section~\ref{sec:num}. Finally, the main conclusions and perspectives are outlined in Section~\ref{sec:disc}.


\section{Symmetry analysis}
\label{sec:sym}

In general, a linear PDE admits an infinity of conservation laws, with integrating multipliers being solutions to the adjoint PDE \cite{Bluman2010}. Here we provide an interesting conservation law, which was found using the \textsc{GeM} \textsc{Maple} package \cite{Cheviakov2007}:
\begin{equation*}
  \Bigl(\E_{\,1}\,\bigl(-\frac{\gamma\,x}{\eta}\bigr)\,p\Bigr)_{\,t}\ +\ \G_{\,x}\ =\ 0\,,
\end{equation*}
where $\E_{\,1}\,(z)\ \eqdef\ \int_{\,1}^{\,+\infty}\,\frac{\ue^{-\,t\,z}}{t}\;\ud t$ is the so-called \emph{exponential integral} function \cite{Abramowitz1965} and the flux $\G$ is defined as
\begin{equation*}
  \G\,(x,\,p)\ \eqdef\ -\eta\,\ue^{\,\frac{\gamma\,x}{\eta}}\,p\ -\ x\,\Bigl(\gamma\,\E_{\,1}\,\bigl(-\frac{\gamma\,x}{\eta}\bigr)\,p\ +\ \eta\,\E_{\,1}\,\bigl(-\frac{\gamma\,x}{\eta}\bigr)\,p_{\,x}\Bigr)\,.
\end{equation*}

The symmetry group of point transformations can be computed using \textsc{GeM} package as well. The infinitesimal generators are given below:
\begin{align*}
  \X_{\,1}\ &=\ \D_{\,t}\,, \\
  \X_{\,2}\ &=\ p\,\D_{\,p}\,, \\
  \X_{\,3}\ &=\ -\,\frac{\ue^{\,-\,\gamma\,t}}{\gamma}\;\D_{\,t}\ +\ \ue^{\,-\,\gamma\,t}\,x\,\D_{\,x}\ -\ \ue^{\,-\,\gamma\,t}\,p\,\D_{\,p}\,, \\
  \X_{\,4}\ &=\ \frac{\ue^{\,\gamma\,t}}{\gamma}\;\D_{\,t}\ +\ \ue^{\,\gamma\,t}\,x\,\D_{\,x}\ -\ \frac{\gamma\,\ue^{\,\gamma\,t}}{\eta}\;x\,p\,\D_{\,p}\,, \\
  \X_{\,5}\ &=\ \ue^{\,(\gamma\ +\ c)\,t}\,\M\Bigl(1\ +\ \frac{c}{\gamma},\,1,\,\frac{\gamma\,x}{\eta}\Bigr)\,\ue^{\,-\,\frac{\gamma\,x}{\eta}}\;\D_{\,p}\,, \\
  \X_{\,6}\ &=\ \ue^{\,(\gamma\ +\ c)\,t}\,\U\Bigl(1\ +\ \frac{c}{\gamma},\,1,\,\frac{\gamma\,x}{\eta}\Bigr)\,\ue^{\,-\,\frac{\gamma\,x}{\eta}}\;\D_{\,p}\,,
\end{align*}
where $c\ \in\ \R\,$, $\M(a,\,b,\,z)$ and $\U(a,\,b,\,z)$ are \textsc{Kummer} special functions \cite{Kummer1837, Abramowitz1965} (see also Appendix~\ref{app:kum}). The corresponding point transformations, which map solutions of \eqref{eq:feller} into other solutions, can be readily obtained by integrating several ODE systems (we do not provide integration details here):
\begin{equation*}
  \bigl(\,t,\,x,\,p\,\bigr)\ \mapsto\ \bigl(\,t\ +\ \eps_{\,1},\,x,\,p\,\bigr)\,,
\end{equation*}
\begin{equation*}
  \bigl(\,t,\,x,\,p\,\bigr)\ \mapsto\ \bigl(\,t,\,x,\,\ue^{\,\eps_{\,2}}\,p\,\bigr)\,,
\end{equation*}
\begin{equation*}
  \bigl(\,t,\,x,\,p\,\bigr)\ \mapsto\ \Bigl(\frac{1}{\gamma}\;\ln\bigl(\eps_{\,3}\,\gamma\ +\ \ue^{\,\gamma\,t}\bigr),\,\frac{\ue^{\,\gamma\,t}}{\eps_{\,3}\,\gamma\ +\ \ue^{\,\gamma\,t}}\;x,\,\bigl(1\ +\ \eps_{\,3}\,\gamma\,\ue^{\,-\,\gamma\,t}\bigr)\,p\Bigr)\,,
\end{equation*}
\begin{equation*}
  \bigl(\,t,\,x,\,p\,\bigr)\ \mapsto\ \Bigl(t\ -\ \frac{1}{\gamma}\;\ln\bigl(1\ -\ \eps_{\,4}\,\gamma\,\ue^{\gamma\,t}\bigr)\,, \frac{x}{1\ -\ \eps_{\,4}\,\gamma\,\ue^{\,\gamma\,t}},\,\ue^{-\,\dfrac{\eps_{\,4}\,\gamma^{\,2}\,x\,\ue^{\,\gamma\,t}}{\eta\,(1\ -\ \eps_{\,4}\,\gamma\,\ue^{\,\gamma\,t})}}\cdot p\Bigr)\,,
\end{equation*}
\begin{equation*}
  \bigl(\,t,\,x,\,p\,\bigr)\ \mapsto\ \Bigl(t,\,x,\,p\ +\ \eps_{\,5}\,\M\bigl(1\ +\ \frac{c}{\gamma},\,1,\,\frac{\gamma\,x}{\eta}\bigr)\,\ue^{-\,\frac{\gamma\,x}{\eta}\ +\ (\gamma\ +\ c)\,t}\Bigr)\,,
\end{equation*}
\begin{equation*}
  \bigl(\,t,\,x,\,p\,\bigr)\ \mapsto\ \Bigl(t,\,x,\,p\ +\ \eps_{\,6}\,\U\bigl(1\ +\ \frac{c}{\gamma},\,1,\,\frac{\gamma\,x}{\eta}\bigr)\,\ue^{-\,\frac{\gamma\,x}{\eta}\ +\ (\gamma\ +\ c)\,t}\Bigr)\,.
\end{equation*}
The first symmetry is the time translation. The second one is the scaling of the dependent variable (the governing equation is linear). Symmetries $3$ and $4$ are exponential scalings. Two last symmetries express the fact that we can always add to the solution a particular solution to the homogeneous equation to obtain another solution. For instance, the solutions invariant under time translations ($\X_{\,1}$) are steady states and their general form is the following:
\begin{equation}\label{eq:steady}
  p\,(x)\ =\ \ue^{\,-\,\frac{\gamma\,x}{\eta}}\Bigl(\Co_{\,1}\,\E_{\,1}\,\bigl(-\frac{\gamma\,x}{\eta}\bigr)\ +\ \Co_{\,2}\Bigr)\,,
\end{equation}
where $\Co_{\,1,\,2}$ are `arbitrary' constants, which have to be determined from imposed conditions. Of course, they should be chosen so that the resulting steady solution is a valid probability distribution. It is not difficult to check that the imposed flux $\F$ on the steady state solution is equal to $\Co_{\,1}\,\eta\,$. Please, notice also an important property of the exponential integral function, which is useful in manipulating its values for negative arguments:
\begin{equation*}
  \lim_{\delta\ \to\ 0} \E_{\,1}\,(x\ \pm\ \ui\,\delta)\ =\ \E_{\,1}\,(x)\ \mp\ \ui\,\pi\ \equiv\ -\,\Ei(-x)\ \mp\ \ui\,\pi\,, \qquad
  x\ <\ 0\,,
\end{equation*}
where $\Ei(z)$ is the following exponential integral:
\begin{equation*}
  \Ei(z)\ \eqdef\ -\,\int_{\,-\infty}^{\,z}\,\frac{\ue^{\,t}}{t}\;\ud t\,.
\end{equation*}
Thus, $\Ei(x)\ \equiv\ -\E_{\,1}\,(-x)$ for $x\ <\ 0\,$. The last relation can be also extended to the entire complex plain:
\begin{equation*}
  \Ei(z)\ \equiv\ -\E_{\,1}\,(-z)\ +\ \frac{1}{2}\;\ln z\ -\ \frac{1}{2}\;\ln\Bigl(\,\frac{1}{z}\,\Bigr)\ -\ \ln(-z)\,, \qquad z\ \in\ \C\,.
\end{equation*}
We provide here also the general solutions invariant under the symmetry $(\X_{\,3})\,$:
\begin{equation*}
  p\,(x,\,t)\ =\ \Bigl(\Co_{\,2}\ -\ \Co_{\,1}\,t\ +\ \frac{\Co_{\,1}}{\gamma}\;\ln x\Bigr)\,\ue^{-\,\frac{\gamma\,x}{\eta}}
\end{equation*}
and under symmetry $(\X_{\,4})\,$:
\begin{equation*}
  p\,(x,\,t)\ =\ \Bigl(\Co_{\,1}\ +\ \Co_{\,2}\,t\ +\ \frac{\Co_{\,2}}{\gamma}\;\ln x\Bigr)\,\ue^{\,\gamma\,t}\,.
\end{equation*}
These solutions might be used, for example, to validate numerical codes.


\section{Reformulation}
\label{sec:reform}

By following the lines of \cite[Chapter~7]{Gosse2013}, we are going to rewrite \textsc{Feller}'s Equation~\eqref{eq:feller} with the so-called \textsc{Lagrangian} or \emph{material} variables. The main advantage of this formulation consists in the fact that we can handle infinite domains \emph{without} any truncations, transformations, \etc It becomes possible to carry computations in infinite domains. Our domain is semi-infinite ($x\ \in\ \R^{\,+}$) with the left boundary $x\ =\ 0$ being a reflection point.

As the first step, we introduce the distribution function associated to the probability density $p\,(x,\,t)\,$:
\begin{equation}\label{eq:78}
  \P\,(x,\,t)\ \eqdef\ \int_{\,0}^{\,x}\,p\,(\xi,\,t)\,\ud\xi\,.
\end{equation}
The same can be done for the initial condition as well:
\begin{equation*}
  \P_{\,0}(x)\ \eqdef\ \int_{\,0}^{\,x}\,p_{\,0}\,(\xi)\,\ud\xi\,, \qquad
  p_{\,0}\ \in\ W_{\,loc}^{\,1,\,1}\,(\R^{\,+})\,.
\end{equation*}
We notice also two obvious properties of the function $\P\,(x,\,t)\,$:
\begin{equation*}
  \partial_{\,x}\P\,(x,\,t)\ \equiv\ p\,(x,\,t)\,, \qquad
  \lim_{x\ \to\ 0} \P\,(x,\,t)\ =\ 0\,, \qquad
  \lim_{x\ \to\ +\infty} \P\,(x,\,t)\ =\ 1\,.
\end{equation*}
Due to the positivity preservation property (see Appendix~\ref{app:pos}) the function $\P\,(x,\,t)$ is non-decreasing in variable $x\,$. Thus, we can define its pseudo-inverse\footnote{This mapping is sometimes called the \emph{reciprocal mapping} \cite{Gosse2013} or an \emph{order preserving string} \cite{Brenier2004}.}:
\begin{equation*}
  \x\,:\ [\,0,\,1\,]\times\R^{\,+}\ \mapsto\ \R^{\,+}\,,
\end{equation*}
which is defined as
\begin{equation*}
  \x\,(\bar{\P},\,t)\ \eqdef\ \inf\;\{\,\xi\ \in\ \R^{\,+}\ \vert\ \P\,(\xi,\,t)\ =\ \bar{\P}\,\}\,.
\end{equation*}
Similarly, the initial condition does possess a pseudo-inverse as well:
\begin{equation}\label{eq:initp}
  \x_{\,0}\,(\bar{\P})\ \eqdef\ \inf\;\{\,\xi\ \in\ \R^{\,+}\ \vert\ \P_{\,0}\,(\xi)\ =\ \bar{\P}\,\}\,,
\end{equation}
such that $\x\,(\bar{\P},\,0)\ \equiv\ \x_{\,0}\,(\bar{\P})\,$.

If \textsc{Feller} equation \eqref{eq:feller} holds in the sense of distributions, then the following equation holds as well:
\begin{equation}\label{eq:shit}
  \P_{\,t}\ -\ x\,\Bigl[\,\gamma\,\P\ +\ \eta\,\P_{\,x}\,\Bigr]_{\,x}\ =\ 0\,.
\end{equation}
along with the initial condition
\begin{equation*}
  \P\,(x,\,0)\ =\ \P_{\,0}(x)\,.
\end{equation*}
Equation \eqref{eq:shit} can be readily obtained by exploiting the obvious property $p\,(x,\,t)\ =\ \partial_{\,x}\,\P\,(x,\,t)\,$. In Appendix~\ref{app:pos} we show that zero value of the solution $p\,(x,\,t)$ is repulsive. Thus, $\partial_{\,x}\,\P\,(x,\,t)\ \equiv\ p\,(x,\,t)\ >\ 0\,$, $\ \forall\, (x,\,t)\ \in\ \bigl(\R^{\,+}\bigr)^{\,2}\,$. Thus, the implicit function theorem \cite{Zorich2004, zorich} guarantees the existence of derivatives of the inverse mapping $\x\,(\bar{\P},\,t)\,$. Let us compute them by differentiating with respect to $\bar{\P}$ and $t$ the following obvious identity:
\begin{equation*}
  \P\,(\x\,(\bar{\P},\,t),\,t)\ \equiv\ \bar{\P}\,.
\end{equation*}
Thus, one can easily show that
\begin{equation*}
  \pd{\x}{\P}\ =\ \frac{1}{\partial_{\,x}\,\P}\,, \qquad
  \pd{\x}{t}\ =\ -\,\frac{\partial_{\,t}\,\P}{\partial_{\,x}\,\P}\,.
\end{equation*}
Using these expressions of partial derivatives, we derive the following evolution equation for the inverse mapping $\x\,(\P,\,\cdot)\,$:
\begin{equation}\label{eq:lagr}
  \bigl(\ue^{\,\gamma\,t}\,\x\bigr)_{\,t}\ +\ \x\cdot\Bigl[\,\eta\,\ue^{\,\gamma\,t}\,\Bigl(\pd{\x}{\P}\Bigr)^{\,-1}\,\Bigr]_{\,\P}\ =\ 0\,.
\end{equation}
The last equation can be rewritten also by introducing a new dynamic variable:
\begin{equation}\label{eq:xy}
  \y\,(\P,\,t)\ \eqdef\ \ue^{\,\gamma\,t}\,\x\,(\P,\,t)\,, \qquad
  \y\,(\P,\,0)\ \equiv\ \x\,(\P,\,0)\,.
\end{equation}
It is not difficult to see that Equation~\eqref{eq:lagr} becomes:
\begin{equation}\label{eq:y}
  \y_{\,t}\ +\ \y\cdot\Bigl[\,\eta\,\ue^{\,\gamma\,t}\,\Bigl(\pd{\y}{\P}\Bigr)^{\,-1}\,\Bigr]_{\,\P}\ =\ 0\,.
\end{equation}
The last equation will be solved numerically in the following Section.


\subsection{Numerical discretization}

Earlier we derived Equation~\eqref{eq:lagr}, which governs the dynamics of the pseudo-inverse mapping $\x\,(\P,\,\cdot)\,$. The initial condition for Equation~\eqref{eq:lagr} is given by the pseudo-inverse \eqref{eq:initp} of the initial condition $\P_{\,0}\,(x)\,$. We discretize Equation~\eqref{eq:lagr} with an explicit discretization in time since it yields the most straightforward implementation.

The first step in our algorithm consists in choosing the initial sampling interval. We make this choice depending on the provided initial condition. Typically, we want to sample only where it is needed. Thus, it seems reasonable to choose the \emph{initial} segment $[\,0,\,\ell_{\,0}\,]$ with $\ell_{\,0}$ being the leftmost location such that
\begin{equation*}
  1\ -\ \P_{\,0}\,(\ell_{\,0})\ <\ \mathsf{tol}\,.
\end{equation*}
In simulations presented below we chose $\mathsf{tol}\ \sim\ \O\,(10^{\,-5})\,$. Then, we choose the initial sampling $\bigl\{\x_{\,k}^{\,0}\bigr\}_{\,k\,=\,0}^{\,N}\ \in\ [\,0,\,\ell_{\,0}\,]\ \subseteq\ \R^{\,+}\,$, with $\x_{\,0}^{\,0}\ =\ 0$ and $\x_{\,N}^{\,0}\ =\ \ell_{\,0}\,$. It is desirable that the initial sampling be adapted to the initial condition, since errors made initially cannot be corrected later. We define also $\P_{\,k}\ =\ \P_{\,0}\,(\x_{\,k}^{\,0})\,$. We stress out that $\bigl\{\P_{\,k}\bigr\}_{\,k\,=\,0}^{\,N}$ stand for a discrete cumulative mass variable and, thus, they are time independent.

More generally, we introduce the following notation:
\begin{equation*}
  \x_{\,k}^{\,n}\ \eqdef\ \x\,(\P_{\,k},\,t^{\,n})\,, \qquad k\ =\ 0,\,1,\,\ldots,\,N\,,
\end{equation*}
with $t^{\,n}\ \eqdef\ n\,\Delta t\,$, $n\ \in\ \N$ and $\Delta t\ >\ 0$ is a chosen time step\footnote{We present our algorithm with a constant time step for the sake of simplicity. However, in realistic simulations presented in Section~\ref{sec:num} the time step will be chosen adaptively and automatically to meet the stability and accuracy requirements prescribed by the user.}. We introduce also similar notation for the dynamic variable:
\begin{equation*}
  \y_{\,k}^{\,n}\ \eqdef\ \y\,(\P_{\,k},\,t^{\,n})\,, \qquad \y_{\,k}^{\,0}\ \equiv\ \x_{\,k}^{\,0}\,, \qquad k\ =\ 0,\,1,\,\ldots,\,N\,.
\end{equation*}
Now we can state the fully discrete scheme for Equation~\eqref{eq:y}:
\begin{equation}\label{eq:fd}
  \frac{\y_{\,k}^{\,n\,+\,1}\ -\ \y_{\,k}^{\,n}}{\Delta t}\ +\ \eta\,\ue^{\,\gamma\,t^{\,n}}\;\frac{\y_{\,k}^{\,n}}{\Delta \P_{\,k}}\;\biggl\{\frac{\Delta\P_{k\,+\,\half}}{\y_{\,k\,+\,1}^{\,n}\ -\ \y_{\,k}^{\,n}}\ -\ \frac{\Delta\P_{k\,-\,\half}}{\y_{\,k}^{\,n}\ -\ \y_{\,k\,-\,1}^{\,n}}\biggr\}\,,
\end{equation}
with $n\ \geq\ 0\,$, $k\ =\ 0,\,1,\,\ldots,\,N\,-\,1$ and
\begin{equation*}
  \Delta\P_{k\,+\,\half}\ \eqdef\ \P_{\,k\,+\,1}\ -\ \P_{\,k}\,, \qquad
  \Delta\P_{k\,-\,\half}\ \eqdef\ \P_{\,k}\ -\ \P_{\,k\,-\,1}\,.
\end{equation*}
The quantity $\Delta \P_{\,k}$ can be defined as the arithmetic or geometric mean of two neighbouring discretization steps $\Delta\P_{k\,\pm\,\half}\,$:
\begin{equation*}
  \Delta \P_{\,k}\ \eqdef\ \frac{\Delta\P_{k\,+\,\half}\ +\ \Delta\P_{k\,-\,\half}}{2}\,, \qquad
  \Delta \P_{\,k}\ \eqdef\ \sqrt{\Delta\P_{k\,+\,\half}\cdot\Delta\P_{k\,-\,\half}}\,.
\end{equation*}
To be specific, in our code we implemented the arithmetic mean. The fully discrete scheme can be easily rewritten under the form of a discrete dynamical system:
\begin{equation*}
  \y_{\,k}^{\,n\,+\,1}\ =\ \y_{\,k}^{\,n}\ -\ \eta\,\Delta t\,\ue^{\,\gamma\,t^{\,n}}\;\frac{\y_{\,k}^{\,n}}{\Delta \P_{\,k}}\;\biggl\{\frac{\Delta\P_{k\,+\,\half}}{\y_{\,k\,+\,1}^{\,n}\ -\ \y_{\,k}^{\,n}}\ -\ \frac{\Delta\P_{k\,-\,\half}}{\y_{\,k}^{\,n}\ -\ \y_{\,k\,-\,1}^{\,n}}\biggr\}\,, \quad n\ \geq\ 0\,.
\end{equation*}

\begin{remark}
We would like to say a few words about the implementation of boundary conditions. First of all, no boundary condition is required on the left side, where $\x_{\,0}^{\,n}\ =\ \y_{\,0}^{\,n}\ \equiv\ 0\,$. On the right boundary we prefer to impose the homogeneous \textsc{Neumann}-type boundary condition, which yields the exact `mass' conservation at the discrete level as well. Namely, at the rightmost cell we have the following fully discrete scheme:
\begin{equation*}
  \y_{\,N}^{\,n\,+\,1}\ =\ \y_{\,N}^{\,n}\ +\ \eta\,\Delta t\,\ue^{\,\gamma\,t^{\,n}}\;\frac{\y_{\,N}^{\,n}}{\Delta \P_{\,N}}\cdot\frac{\Delta\P_{N\,-\,\half}}{\y_{\,N}^{\,n}\ -\ \y_{\,N\,-\,1}^{\,n}}\,, \qquad n\ \geq\ 0\,,
\end{equation*}
with $\Delta \P_{\,N}\ \eqdef\ \dfrac{\P_{\,N}\ -\ \P_{\,N\,-\,2}}{2}\,$. As a result, we obtain the exact conservation of `mass' at the discrete level:
\begin{equation*}
  \sum_{k\,=\,0}^{\,N}\,\Delta \P_{\,k}\,\x_{\,k}^{\,n}\ \equiv\ \sum_{k\,=\,0}^{\,N}\,\Delta \P_{\,k}\,\x_{\,k}^{\,0}\,, \qquad \forall n\ \in\ \N\,.
\end{equation*}
\end{remark}

To summarize, our numerical strategy consists in:
\begin{enumerate}
  \item We compute the pseudo-inverse of the initial data $p_{\,0}\,(x)$ to obtain $\x\,(\P,\,0)\ \equiv\ \y\,(\P,\,0)\,$.
  \item This initial condition $\y\,(\P,\,0)$ is evolved in (discrete) time using an \emph{explicit} marching scheme in order to obtain numerical approximation to $\y\,(\P,\,t)\,$, $t\ >\ 0\,$.
  \item The variable $\x\,(\P,\,t)$ is recovered by inverting \eqref{eq:xy}, \ie
  \begin{equation*}
    \x\,(\P,\,t)\ =\ \ue^{\,-\gamma\,t}\,\y\,(\P,\,t)\,.
  \end{equation*}
  \item Thanks to \eqref{eq:78} we can deduce the values of $p\,\bigl(\x(\P,\,t),\,t\bigr)\ \in\ [\,0,\,1\,]\,$ by applying a favourite finite difference formula\footnote{In our code we employed the simplest forward finite differences and it lead satisfactory results. This point can be easily improved when necessary.}.
\end{enumerate}
Working with the pseudo-inverse allows to overcome the issue of the retention phenomenon, which manifests as the expanding support of $p\,(x,\,t)$ for positive (and possibly large) times $t\ >\ 0\,$, $t\ \gg\ 1\,$, since the computational domain was transformed to $[\,0,\,1\,]\,$. This method is the \textsc{Lagrangian} counterpart of the moving mesh technique in the \textsc{Eulerian} setting \cite{Khakimzyanov2015a, Khakimzyanov2015b}.

A simple \textsc{Matlab} code, which implements the scheme we described hereinabove, is freely available for reader's convenience at this URL address:
\smallskip

\url{https://github.com/dutykh/Feller/}


\section{Numerical results}
\label{sec:num}

In this Section we validate and illustrate the application of the proposed algorithm on several examples. However, first we begin with a straightforward validation test. The only difference with the proposed above algorithm is that we are using a higher order adaptive time stepping for our practical simulations. The explicit first order scheme was used to simplify the presentation. In practice, much more sophisticated time steppers can be used. For instance, we shall employ the explicit embedded \textsc{Dormand}--\textsc{Prince} \textsc{Runge}--\textsc{Kutta} pair $(4,\,5)$ \cite{Dormand1980} implemented in \textsc{Matlab} under the \texttt{ode45} routine \cite{Shampine1997}. Conceptually, this method is similar to the explicit \textsc{Euler} scheme presented above. It conserved all good properties we showed, but it provides additionally the higher accuracy order and totally automatic adaptivity of the time step, which matches very well with adaptive features of the \textsc{Lagrangian} scheme in space. The values of absolute and relative tolerances used in the time step choice are systematically reported below.


\subsection{Steady state preservation}
\label{sec:wb}

\begin{figure}
  \centering
  \includegraphics[width=0.59\textwidth]{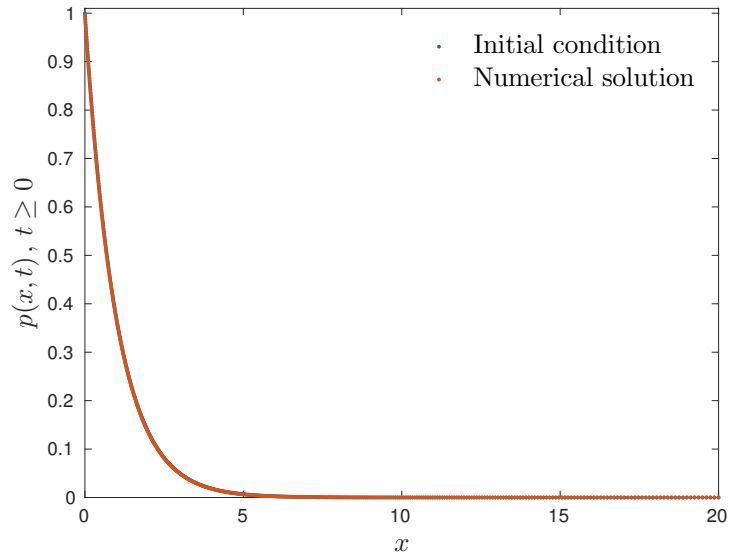}
  \caption{\small\em Comparison of a steady state solution of class \eqref{eq:steady} at $t\ =\ 0$ and at $t\ =\ T\ =\ 10\,$. They are indistinguishable up to the graphical resolution, which validates the solver.}
  \label{fig:steady}
\end{figure}

In order to validate the numerical algorithm we have at our disposal a family of steady state solutions \eqref{eq:steady}. Hence, if we take such a solution as an initial condition, normally the algorithm has to keep it intact under the discretized dynamics. The parameters $\eta\,$, $\gamma$ of the equation, those of the steady solution and numerical parameters used in our computation are reported in Table~\ref{tab:params0}. The initial condition at $t\ =\ 0$ along with the final state at $t\ =\ T$ are shown in Figure~\ref{fig:steady}. Up to graphical resolution they coincide completely. We can easily check that during the whole simulation the points barely moved, \ie
\begin{equation*}
  \norm{\x(\cdot,\,T)\ -\ \x(\cdot,\,0)}_{\,\infty}\ \approx\ 0.008577\ldots
\end{equation*}
We can check also other quantities. For instance, $\P\,(\X,\,t)$ is preserved up to the machine precision. If we reconstruct the probability distribution $p\,(x,\,t)\,$, we obtain:
\begin{equation*}
  \norm{p\,(\cdot,\,T)\ -\ p\,(\cdot,\,0)}_{\,\infty}\ \approx\ 0.003051\ldots
\end{equation*}
The last error comes essentially from the fact that we apply simple first order finite difference to reconstruct the variable $p\,(x,\,t)$ from its primitive $\P\,(x,\,t)\,$. We can improve this point, but even this simple reconstruction seems to be consistent with the overall scheme accuracy. Thus, this example shows that our implementation of the proposed algorithm is also practically well-balanced \cite{Gosse2013}, since steady state solutions are well preserved.

\begin{table}
  \begin{tabular}{l|c}
  \hline\hline
  \textit{Parameter} & \textit{Value} \\
  \hline\hline
  Drift coefficient, $\gamma$ & $1.0$ \\
  Diffusion coefficient, $\eta$ & $1.0$ \\
  Integration constant, $\Co_{\,1}$ & $0.0$ \\
  Integration constant, $\Co_{\,2}$ & $\dfrac{\eta}{\gamma}\ \equiv\ 1.0$ \\
  Final simulation time, $T$ & $10.0$ \\
  Number of discretization points, $N$ & $500$ \\
  Absolute tolerance, $\mathsf{tol}_{\,\mathrm{a}}$ & $10^{\,-5}$ \\
  Relative tolerance, $\mathsf{tol}_{\,\mathrm{r}}$ & $10^{\,-5}$ \\
  \hline\hline
  \end{tabular}
  \bigskip
  \caption{\small\em Numerical parameters used in the steady state computations.}
  \label{tab:params0}
\end{table}

\subsection{Transient computations}

In this Section we present a couple of extra truly unsteady computations in order to illustrate the capabilities of our method. Namely, we shall simulate the probability distributions emerging from a family of initial conditions (normalized to have the probability distribution):
\begin{equation*}
  p_{\,0}\,(x)\ =\ \frac{\ue^{\,-\frac{x}{\sigma_{\,1}}}\ +\ \ue^{\,-\frac{x\ -\ x_{\,0}}{\sigma_{\,2}}}}{\sigma_{\,1}\ +\ \sigma_{\,2}\,\ue^{\frac{x_{\,0}}{\sigma_{\,2}}}}\,.
\end{equation*}
The primitive of the last distribution can be easily computed as well:
\begin{equation*}
  \P\,(x)\ =\ 1\ -\ \frac{\sigma_{\,1}\,\ue^{\,-\frac{x}{\sigma_{\,1}}}\ +\ \sigma_{\,2}\,\ue^{\,-\frac{x\ -\ x_{\,0}}{\sigma_{\,2}}}}{\sigma_{\,1}\ +\ \sigma_{\,2}\,\ue^{\frac{x_{\,0}}{\sigma_{\,2}}}}\,.
\end{equation*}
We design two different experiments \emph{in silico} to show two completely different behaviours of solutions to \textsc{Feller} equation \eqref{eq:feller} depending on the sign of the drift coefficient $\gamma\,$. These will constitute additional tests for the proposed numerical method. In both cases, the initial positions of particles are chosen according to the logarithmic distribution (\texttt{logspace} function in \textsc{Matlab}) on the segment $[\,0,\,20\,]\,$. This choice is made to represent more accurately the exponentially decaying initial condition since the errors made in the initial condition cannot be corrected in the dynamics.


\subsubsection{Expanding drift}

\begin{figure}
  \centering
  \subfigure[]{\includegraphics[width=0.48\textwidth]{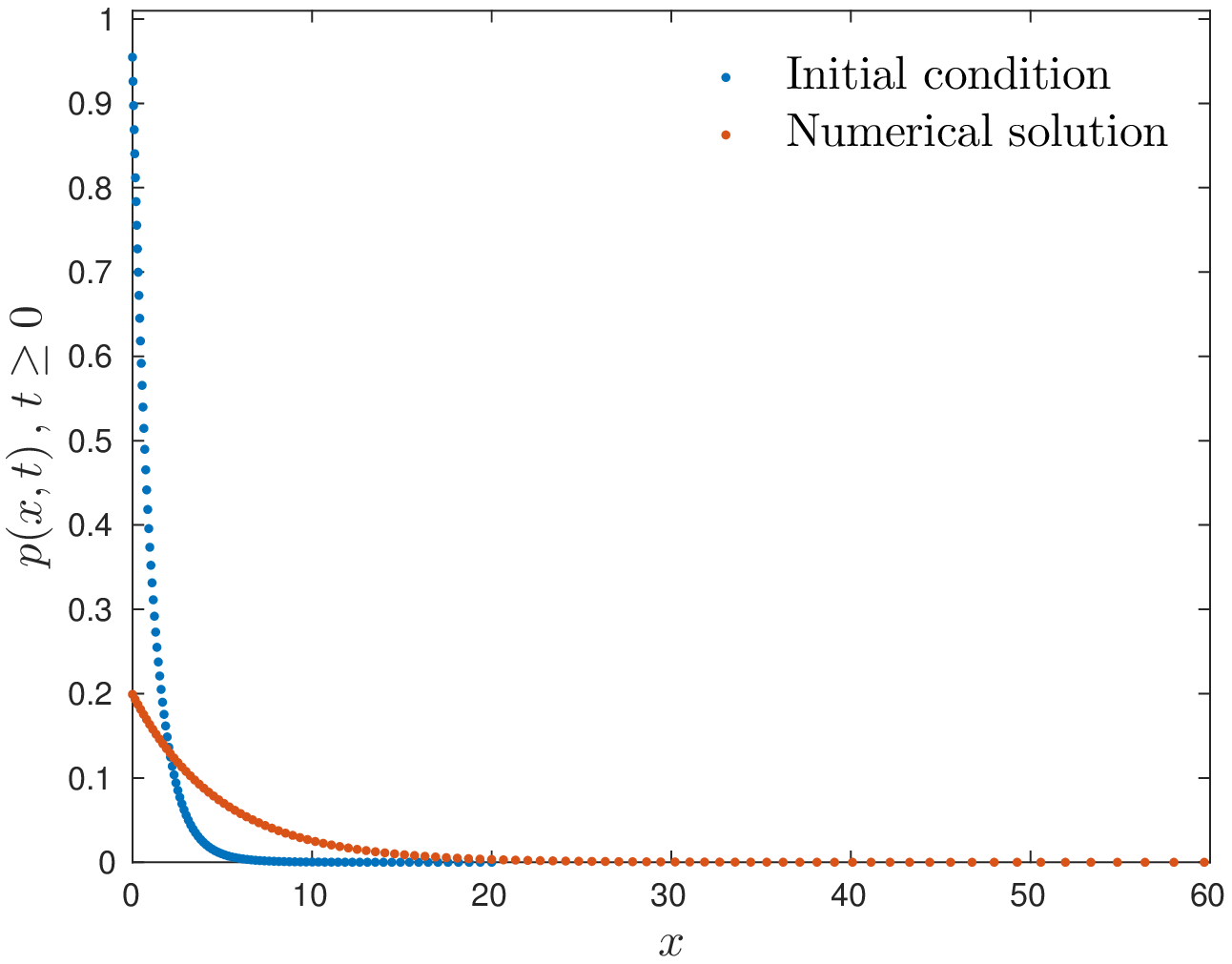}}
  \subfigure[]{\includegraphics[width=0.48\textwidth]{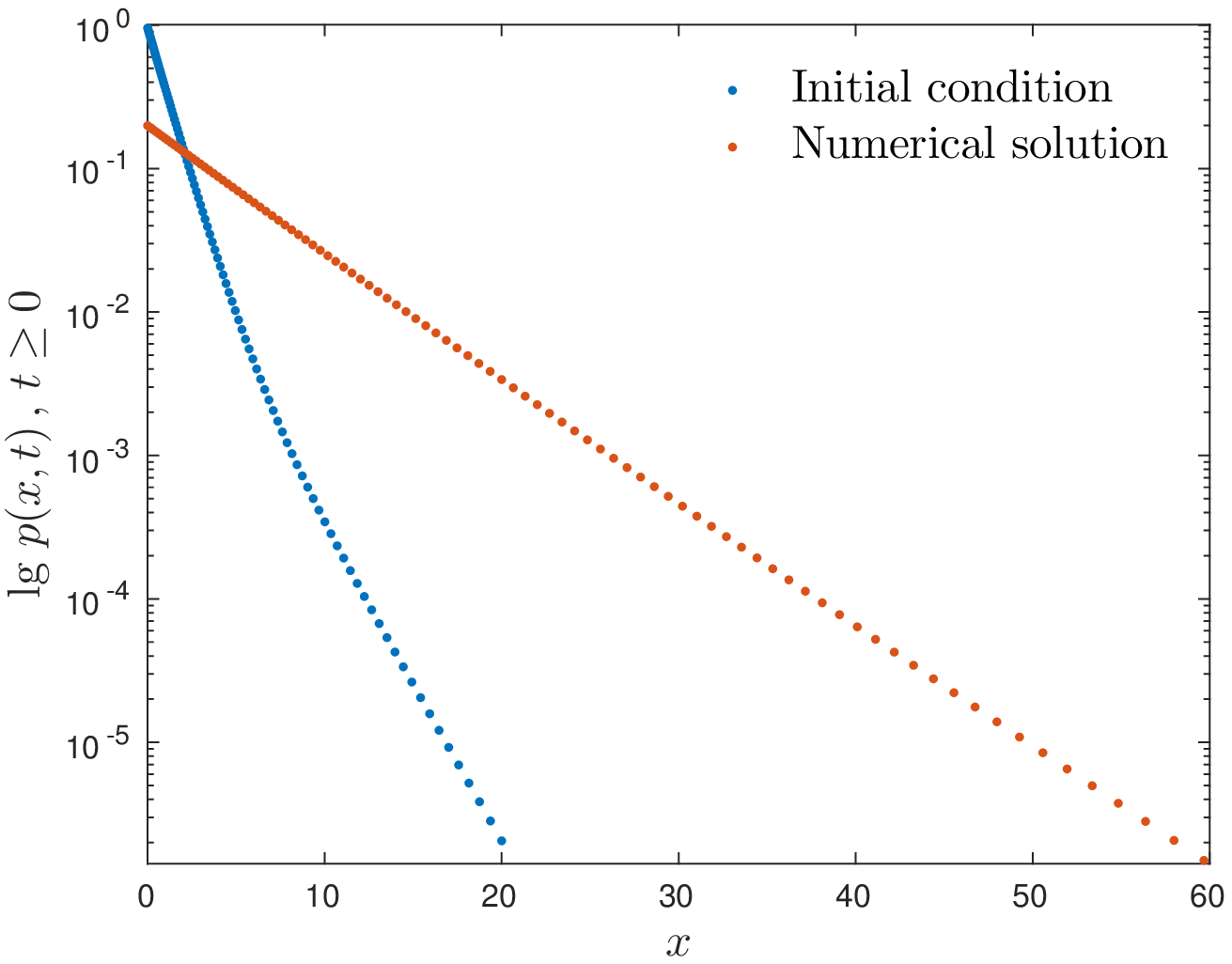}}
  \subfigure[]{\includegraphics[width=0.65\textwidth]{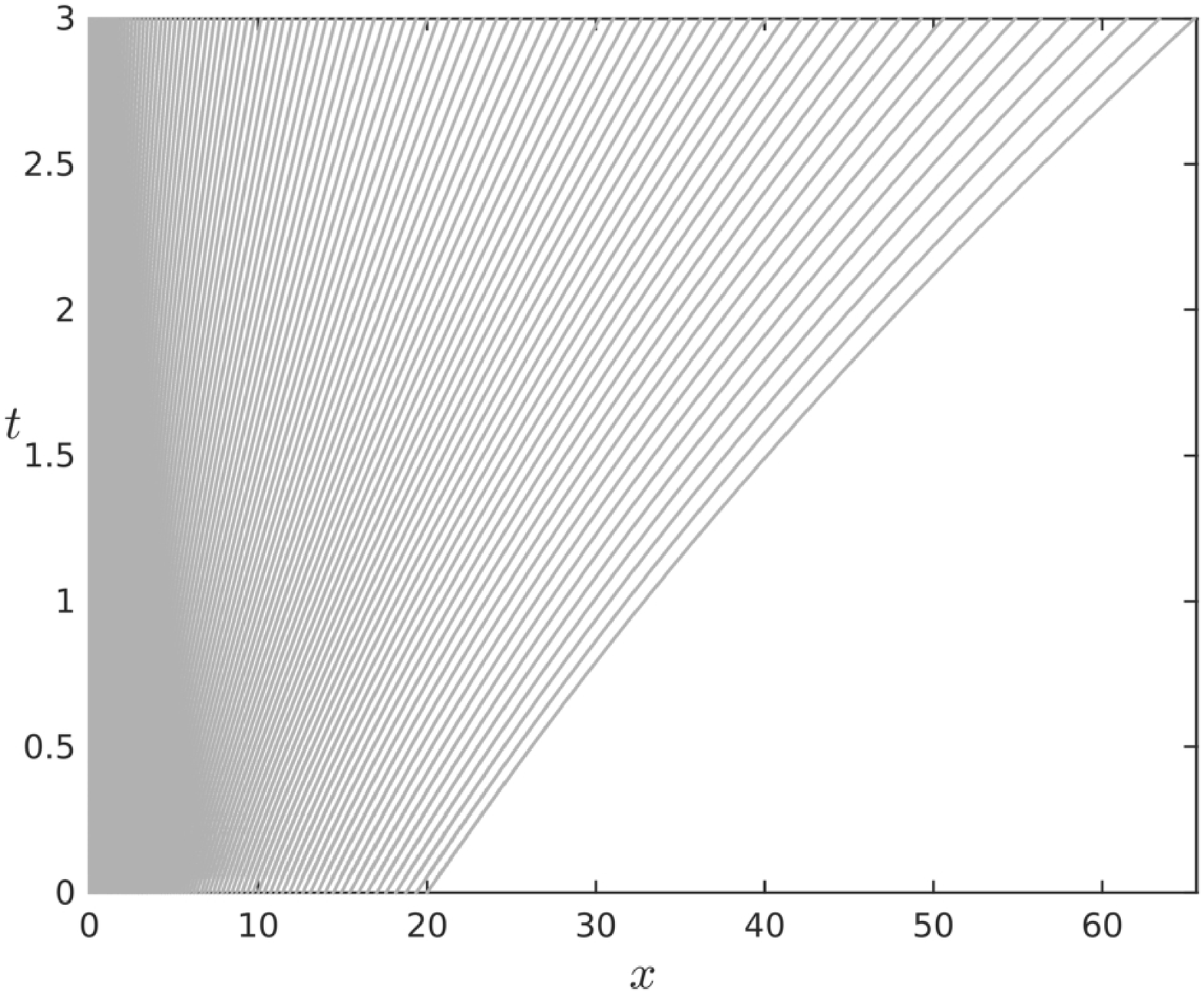}}
  \caption{\small\em Numerical result of the expanding experiment with negative drift $\gamma\ =\ -0.1\ <\ 0\,$: (a) initial and terminal states of the numerical discretized solution; (b) initial and terminal states of the numerical discretized solution in semi-logarithmic coordinates; (c) trajectories of grid nodes. All numerical parameters for this computation are reported in Table~\ref{tab:params1} (middle column).}
  \label{fig:expa}
\end{figure}

If the drift coefficient $\gamma\ <\ 0\,$, this term will cause \textsc{Feller}'s equation \eqref{eq:feller} to transport information in the rightward direction. This situation is quite uncomfortable from the numerical point of view since the initial condition expands quickly towards the (positive) infinity. We submitted our method to this case. All numerical and physical parameters are provided in Table~\ref{tab:params1} (the middle column). The initial condition along with the probability distribution at the end of our simulation are shown simultaneously in Figure~\ref{fig:expa}(\textit{a,b}) in \textsc{Cartesian} and semi-logarithmic coordinates correspondingly. As expected, the support of the initial condition more than triples from $t\ =\ 0$ to $t\ =\ T\ =\ 3.0\,$. We remind that the diffusion and drift coefficients are proportional to $x$ and the scheme handles the growth of these terms automatically. The smooth decay of the solution on the semi-logarithmic plot (see Figure~\ref{fig:expa}(\textit{b})) shows the absence of any numerical instabilities. The trajectory of grid nodes is shown in Figure~\ref{fig:expa}(\textit{c}). One can see that points follow the expansion of the solution. Nevertheless, they concentrate in the areas where the probability takes significant values.

\begin{table}
  \begin{tabular}{l|c|c}
  \hline\hline
  \textit{Parameter} & \textit{Expanding experiment} & \textit{Confining experiment} \\
  \hline\hline
  Drift coefficient, $\gamma$ & $-0.1$ & $0.5$ \\
  Diffusion coefficient, $\eta$ & $1.0$ & $1.0$ \\
  Final simulation time, $T$ & $3.0$ & $12.0$ \\
  Number of discretization points, $N$ & $100$ & $100$ \\
  Initial condition parameter, $\sigma_{\,1}$ & $2.0$ & $2.0$ \\
  Initial condition parameter, $\sigma_{\,2}$ & $1.0$ & $1.0$ \\
  Initial condition parameter, $x_{\,0}$ & $3.0$ & $3.0$ \\
  Absolute tolerance, $\mathsf{tol}_{\,\mathrm{a}}$ & $10^{\,-5}$ & $10^{\,-5}$ \\
  Relative tolerance, $\mathsf{tol}_{\,\mathrm{r}}$ & $10^{\,-5}$ & $10^{\,-5}$ \\
  \hline\hline
  \end{tabular}
  \bigskip
  \caption{\small\em Numerical parameters used in unsteady computations.}
  \label{tab:params1}
\end{table}

\subsubsection{Confining drift}

\begin{figure}
  \centering
  \subfigure[]{\includegraphics[width=0.47\textwidth]{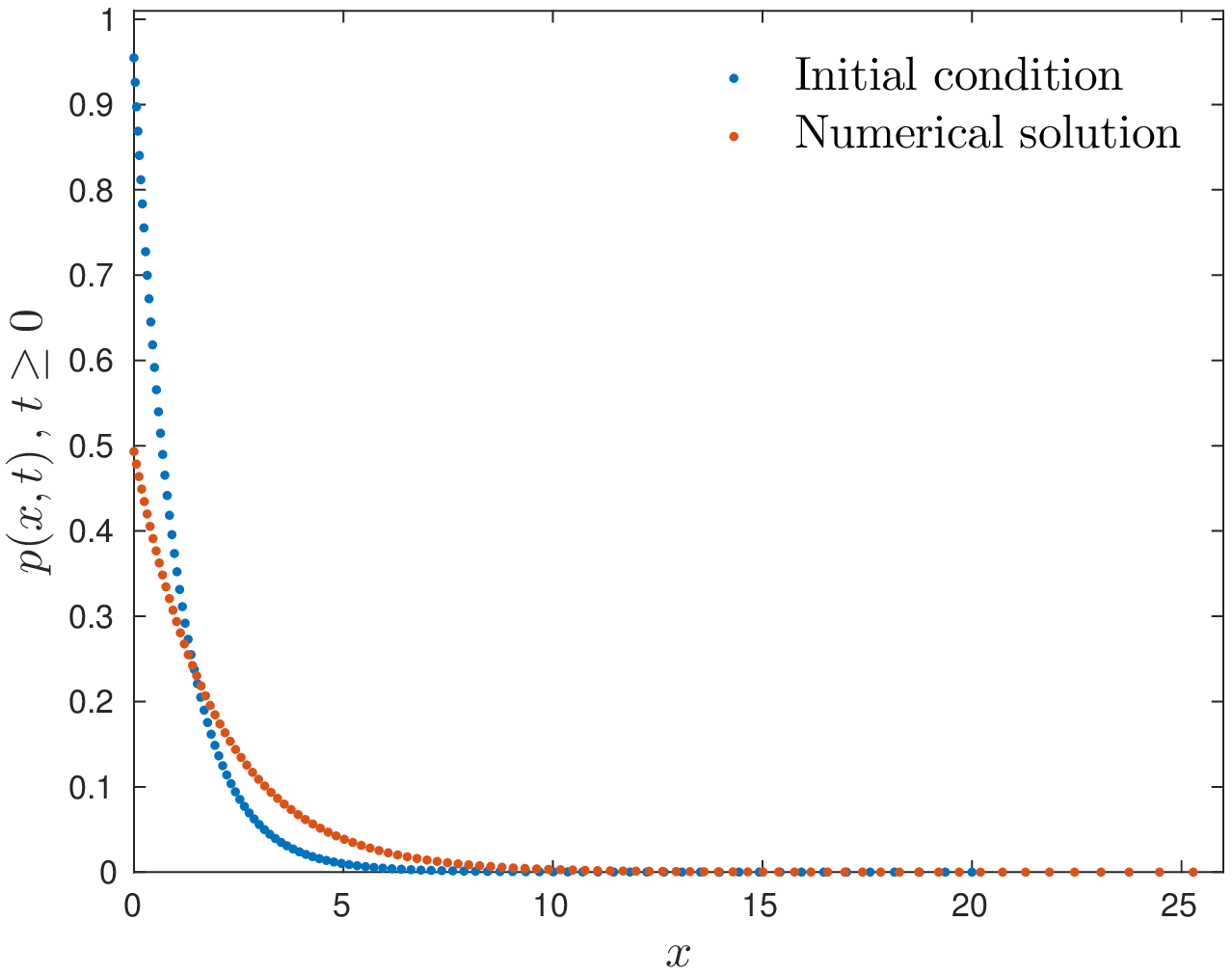}}
  \subfigure[]{\includegraphics[width=0.49\textwidth]{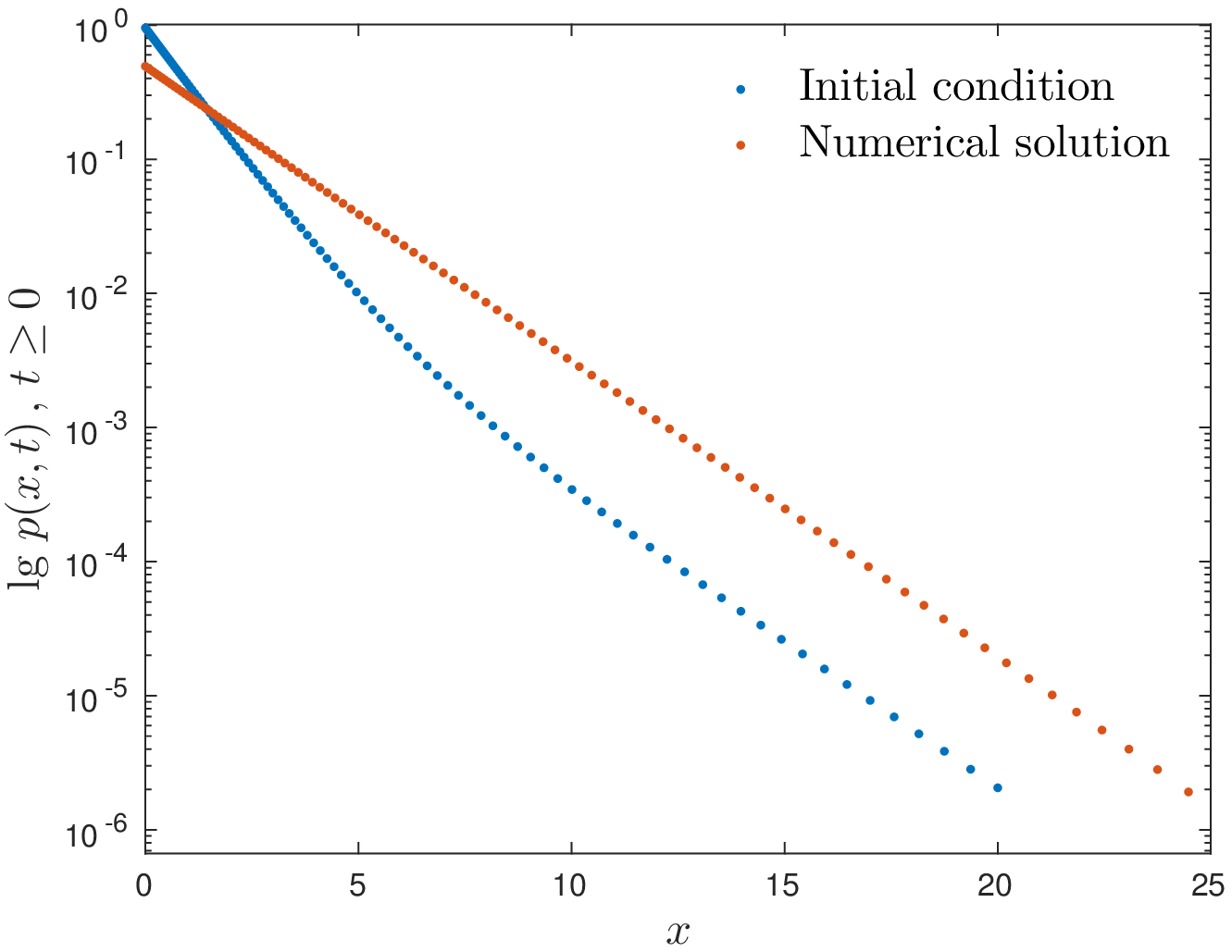}}
  \subfigure[]{\includegraphics[width=0.65\textwidth]{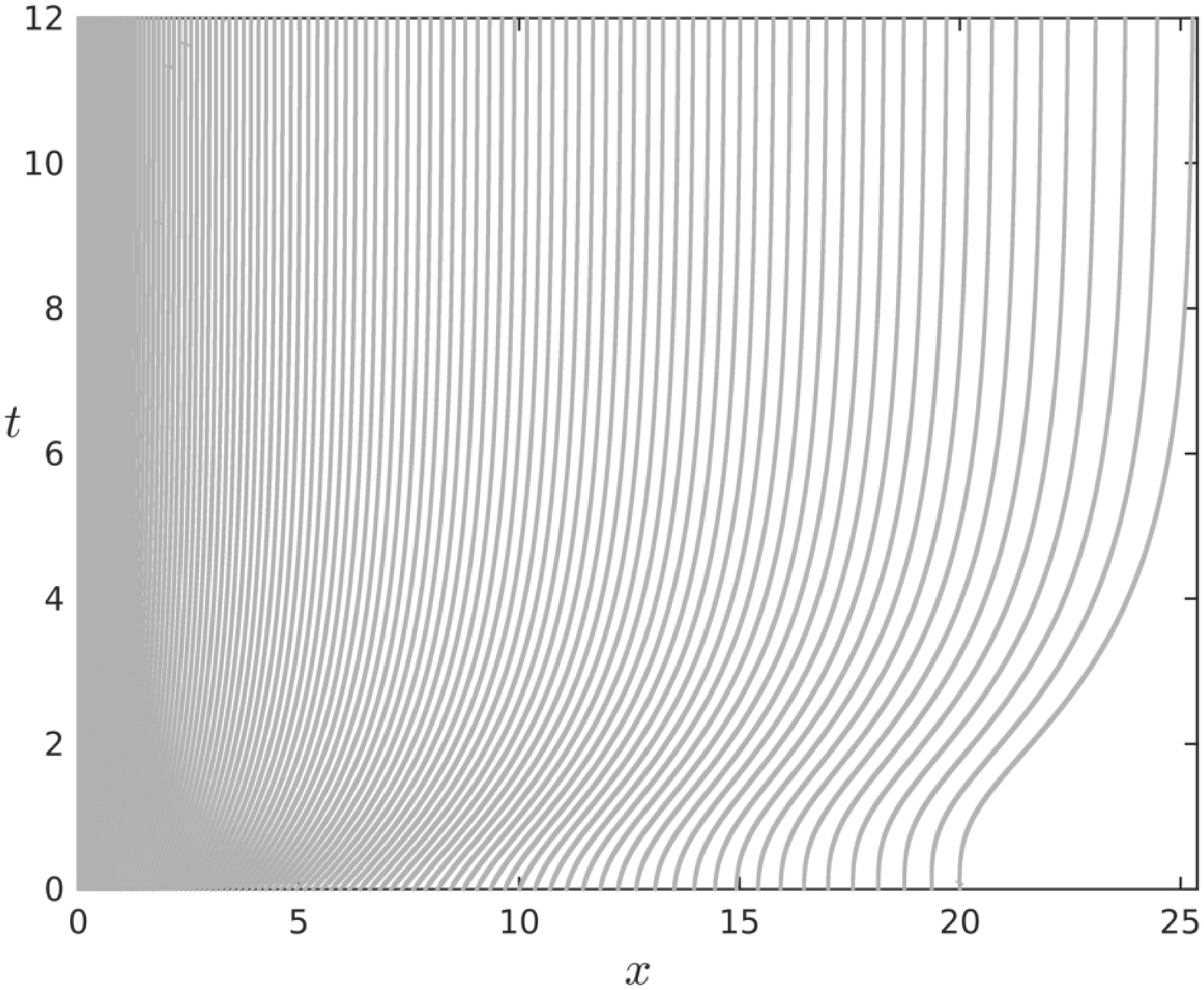}}
  \caption{\small\em Numerical result of the confining experiment with positive drift $\gamma\ =\ 0.5\ >\ 0\,$: (a) initial and terminal states of the numerical discretized solution; (b) initial and terminal states of the numerical discretized solution in semi-logarithmic coordinates; (c) trajectories of grid nodes. All numerical parameters for this computation are reported in Table~\ref{tab:params1} (right column).}
  \label{fig:std}
\end{figure}

In the case of the positive drift coefficient $\gamma\ >\ 0\,$, the \textsc{Feller} dynamics gets even more interesting since we have two competing effects:
\begin{enumerate}
  \item Positive drift moving information towards the origin $x\ =\ 0\,$,
  \item Diffusion spreading solution towards the positive infinity.
\end{enumerate}
It might happen that both effects balance each other and result in a steady solution. Such a simulation is presented in this Section. The numerical and physical parameters are given in Table~\ref{tab:params1} (the right column). The initial condition along with the probability distribution at the end of our simulation are shown simultaneously in Figure~\ref{fig:std}(\textit{a,b}) in \textsc{Cartesian} and semi-logarithmic coordinates respectively. The trajectories of grid nodes are shown in Figure~\ref{fig:expa}(\textit{c}). One can see the rapid initial expansion stage, which is followed by a stabilization process, when we almost achieved the expected steady state. Again, the grid nodes trajectories show excellent adaptive features of the numerical scheme: at the end of the simulation the relative density of nodes turns out to be preserved. The semi-logarithmic plot shown in Figure~\ref{fig:std}(\textit{b}) shows that the numerical solution is free of numerical instabilities. The obtained steady solution will be preserved by discrete dynamics thanks to the well-balanced property demonstrated in Section~\ref{sec:wb}.

\section{Discussion}
\label{sec:disc}

Above we presented some rationale about the discretization, existence and uniqueness theory for \textsc{Feller}'s equation. The main conclusions and perspectives are outlined below.


\subsection{Conclusions}

The celebrated \textsc{Feller} equation was studied mathematically in two seminal papers published by William \textsc{Feller} (1951/1952) in Annals of Mathematics \cite{Feller1951, Feller1952}. In particular, the uniqueness was shown in \cite{Feller1952} without any boundary condition required at $x\ =\ 0\,$. This result is notable and we use it in our study. The main goal of our work was to present an efficient numerical scheme, which were able to handle the unbounded diffusion inherent to Equation~\eqref{eq:feller}. Moreover, the retention phenomenon causes the solution support to expand all the time. Thus, if we wait sufficiently long time, it will reach the computational domain boundaries (since $\R^{\,+}\ \ni\ x$ was truncated to a finite segment to make the simulation \emph{in silico} possible). To overcome this difficulty, we proposed a simple and explicit \textsc{Lagrangian} scheme using the notion of the pseudo-inverse. In this way, we obtained a stable numerical scheme (under not so stringent stability conditions), which turns out to be naturally \emph{adaptive} as well, since particles move together with the growing support (the rightmost particle position defines the support) and they tend to concentrate where it is really needed. The performance of our scheme was illustrated on several examples. We share also the \textsc{Matlab} code so that anybody can check the claims we made hereinabove and use it to solve numerically the \textsc{Feller} equation in their applications:
\smallskip

\url{https://github.com/dutykh/Feller/}


\subsection{Perspectives}

All the numerical schemes and results presented in this paper were in one `spatial' dimension. The \textsc{Feller} equation considered here is $1-$D as well. However, it seems interesting\footnote{Regardless the physical applications.} to consider generalized \textsc{Feller} equations in higher space dimensions and to extend the proposed numerical strategy to the multi-dimensional case as well. Another possible extension direction consists in proposing higher order schemes to capture numerical solutions with better accuracy with the same number of degrees of freedom. On the more theoretical side, we would like to obtain an alternative well-posedness theory for \textsc{Feller} equation by taking a continuous limit in our numerical scheme.


\subsection*{Acknowledgments}
\addcontentsline{toc}{subsection}{Acknowledgments}

The author would like to thank Prof.~Laurent~\textsc{Gosse} (CNR--IAC Rome, Italy) for his invaluable help in understanding his first book. We thank also Prof.~Francesco~\textsc{Fedele} (Georgia Tech, Atlanta, USA) for pointing to us this interesting problem.

\bigskip


\appendix
\section{Positivity preservation}
\label{app:pos}

\begin{figure}
  \centering
  \includegraphics[width=0.89\textwidth]{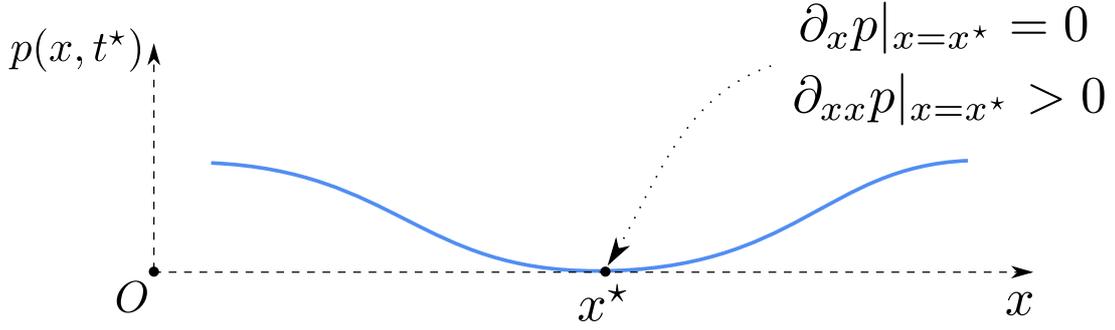}
  \caption{\small\em A schematic representation of the situation where the solution $p\,(x,\,t)$ attains a zero value in some point $x^{\,\star}\ >\ 0\,$.}
  \label{fig:sketch}
\end{figure}

In this Appendix we show the positivity preservation property of the solution $p\,(x,\,t)$ to Equation~\eqref{eq:feller}. We suppose that initially the solution is non-negative, \ie
\begin{equation*}
  p\,(x,\,0)\ =\ p_{\,0}\,(x)\ \geq\ 0\,, \qquad \forall x\ \in\ \R^{\,+}\,.
\end{equation*}
Let us assume that at some positive time $t\ =\ t^{\,\ast}\ >\ 0$ and in some point $x\ =\ x^{\,\ast}\ \in\ \R^{\,+}$ the solution attains zero value, \ie
\begin{equation*}
  p\,(x^{\,\ast},\,t^{\,\ast})\ =\ 0\,.
\end{equation*}
This situation is schematically depicted in Figure~\ref{fig:sketch}.

Equation~\eqref{eq:feller} can be rewritten in the non-conservative form:
\begin{equation*}
  p_{\,t}\ =\ x\,\bigl[\,\gamma\,p_{\,x}\ +\ \eta\,p_{\,x\,x}\,\bigr]\ +\ \gamma\,p\ +\ \eta\,p_{\,x}\,.
\end{equation*}
Taking into account the fact that the point $x^{\,\star}$ is a local minimum ($p_{\,x}\,\vert_{\,x^{\,\star}}\ =\ 0$), where the solution takes zero value ($p\,\vert_{\,x^{\,\star}}\ =\ 0$), the last equation greatly simplifies at this point:
\begin{equation*}
  p_{\,t}\,\vert_{x^{\,\star}}\ =\ x^{\,\star}\,\eta\,p_{\,x\,x}\,\vert_{\,x^{\,\star}}\ >\ 0\,,
\end{equation*}
since in the minimum $p_{\,x\,x}\,\vert_{\,x^{\,\star}}\ >\ 0\,$. Thus, zero values of the solution are \emph{repulsive} and for (at least small) times $t\ >\ t^{\,\star}$ the function $t\ \mapsto\ p\,(x^{\,\star},\,t)$ will be increasing.


\section{`Mass' conservation}
\label{app:norm}

It is straightforward to show that the $L_{\,1}$ norm of the solution to Equations~\eqref{eq:feller} is preserved. Indeed, taking into account that the solution $p\,(x,\,t)\ >\ 0$ is positive for all times $t\ >\ 0$ provided that the initial condition $p_{\,0}\,(x)\ \geq\ 0\,$, $\forall x\ \in\ \R^{\,+}\,$, we have $\abs{p\,(x,\,t)}\ \equiv\ p\,(x,\,t)\,$. By integrating Equation~\eqref{eq:feller}, we have
\begin{equation*}
  \partial_{\,t}\,\int_{\,\R^{\,+}}\,p\,(x,\,t)\,\ud x\ +\ \int_{\,\R^{\,+}}\F_{\,x}\,\ud x\ =\ 0\,.
\end{equation*}
Taking into account that $\F\,\vert_{\,x\,=\,0}\ \equiv\ 0$ and the solution $p\,(x,\,t)$ is decaying sufficiently fast as $x\ \to\ +\,\infty$ together with its derivative, we obtain
\begin{equation*}
  \partial_{\,t}\,\int_{\,\R^{\,+}}\,p\,(x,\,t)\,\ud x\ \equiv\ 0\,.
\end{equation*}
In other words,
\begin{equation*}
  \norm{p\,(\cdot,\,t)}_{\,L_{\,1}}\ =\ \const\,.
\end{equation*}
The last constant can be in general taken equal to one after the appropriate rescaling (provided that the initial condition is non-trivial). It is this scaling, which is assumed throughout the whole text above. This renormalization is consistent with the `physical sense' of the variable $p$ being the density of a probability distribution.


\section{Kummer functions}
\label{app:kum}

The \textsc{Kummer} functions $\M(a,\,b,\,z)$ and $\U(a,\,b,\,z)$ are two linearly independent solutions of the following ordinary differential equation:
\begin{equation*}
  z\,\od{y^{\,2}}{z^{\,2}}\ +\ (b\ -\ z)\;\od{y}{z}\ -\ a\,y\ =\ 0\,.
\end{equation*}
This equation admits two singular points: $z\ =\ 0$ (regular) and $z\ =\ +\infty$ (irregular). There exist a connection between \textsc{Kummer} and hyper-geometric functions \cite{Abramowitz1965}.


\bigskip
\addcontentsline{toc}{section}{References}
\bibliographystyle{abbrv}
\bibliography{biblio}
\bigskip

\end{document}